\documentclass{article}
\usepackage{amssymb}
\usepackage{amsmath}
\usepackage{graphicx}

\newtheorem{definition}{Definition}

\def\d{\partial}
\newcommand{\hf}{\frac12}
\newcommand{\arr}{\rightarrow}

\newcommand{\R}{{\mathbb R}}
\newcommand{\T}{{\mathbb T}}
\newcommand{\Z}{{\mathbb Z}}

\def\fA{\mathfrak{A}}

\newcommand{\dist}{{\operatorname{dist}}}

\begin{document}

\title{On numerical study of attractors of ODEs}
\author{Alexandr Danilin \\ Immanuel Kant Baltic Federal University \\
 ~~~ and ~~~\\
  Iryna Ryzhkova-Gerasymova \\
 V.N. Karazin Kharkiv National University  }
 
 \maketitle

\section{Introduction}
Numerical study of long-time behaviour of solutions of differential equations (ordinary or with partial derivatives) may appear a challenging problem. For example, let's consider the following quasi-hamiltonian 2D system of ODE (Chueshov \cite{Chu99})
\begin{equation}\label{eight_DS}
\left\{
  \begin{array}{l}
    \dot{q}=\frac{\d H}{\d p}-\mu H \frac{\d H}{\d q}, \\
    \dot{p}=-\frac{\d H}{\d q}-\mu H \frac{\d H}{\d p},
  \end{array}
\right.
\end{equation}
with $H(p,q)=\hf p^2+ q^4-q^2$ and $\mu>0$. It generates a dynamical system and it's dynamics is shown on the picture below. It is easy to see, that trajectories of the system cannot cross the separatrix  $\Gamma=\{ (q,p)\, :\;  H(p,q)=0\}$.

\begin{figure}[h]
\begin{center}
\includegraphics[scale=.28]{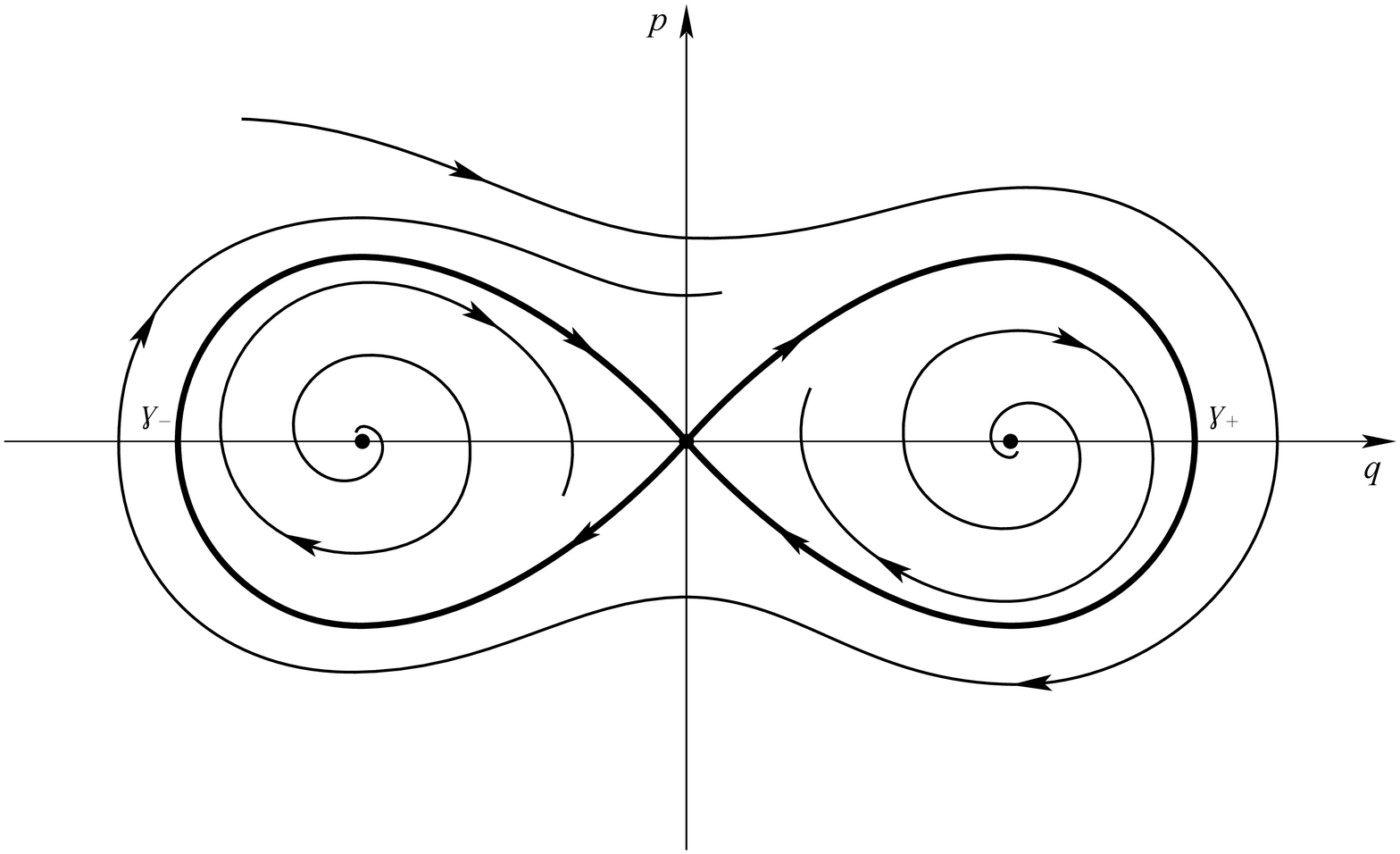}\end{center}
\caption{}
\label{fig:quasi-ham}
\end{figure}

However,  if we simulate numerically an individual trajectory of the system, after some time it may cross the separatirix.

Thus, qualitative behaviour of the numerically simulated solution may be completely different from qualitative behaviour of the real trajectory. Therefor numerical simulations of individual trajectories may give us false idea of long-time behaviour of dynamical system.

In order to overcome this difficulty, the idea is to simulate a bundle of trajectories for a (relatively) short time, rather then individual trajectories for a long time. This idea was presented first by Denlitz etc., and was developed in family of  so-called set-oriented methods for invariant objects. They include algorithms of building covering for stationary points, global attractors, unstable manifolds of stationary points, etc.  We will not give details of these methods here. We give here some examples of covering global attractor of 2D and 3D dynamical systems. Unfortunately, this general method often gives rather rude results. 

For example, we try to construct a cover of global attractor for the system
\begin{eqnarray}
  \dot{x} &=& y, \nonumber \\
  \dot{y} &=& \frac 32 z + \frac{\Gamma-1}{2} x - \frac 12 x^3,  \label{fp_appr}\\
  \dot{z} &=& -y - \frac{11}2 z - \frac{\Gamma-1}{2} x + \frac 12 x^3,  \nonumber
\end{eqnarray}
with $\Gamma=3$.

\begin{figure}[h]
\includegraphics[width=\textwidth]{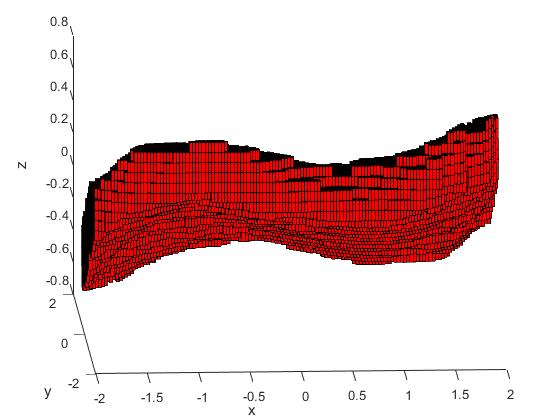}
\caption{}
\end{figure}
In this case subdivision proses was performed 20 times. When we attempted to perform it 30 times, a calculation time and other resources requirements were beyond any acceptable limits.

If we have additional information on the dynamical system, sometimes  we can use other more specific methods form this family  to get better results. For example, system \eqref{fp_appr}  is gradient and it's attractor consists of unstable manifolds, that emanate from stationary points of the system. We can use continuation method from GAIO family in this case an obtain much more accurate result in reasonable time.

\begin{figure}[h]
\includegraphics[width=\textwidth]{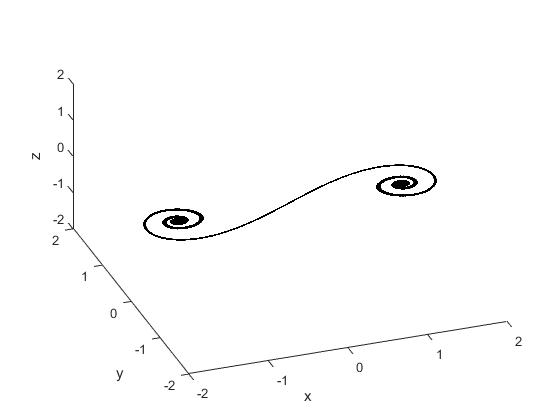}
\caption{}
\end{figure}

In this work the authors give another general method of analysis of asymptotical behaviour of DS, which is based on simulation of bundle of trajectories too. The method is probabilistic in some sense and based on finding domains with high density of trajectories. It leads to approximation of so-called Milnor's attractor and demonstrates high performance. The method is heuristic by now, and we are working on rigorous mathematical justification of it.

The paper is organized as follows. In the 

\section{Definitions and notations}
In this section we give basic definitions of the dynamical systems theory we will use later. More details can be found in \cite{Chu99,Chu2015}.

\begin{definition}\label{de:e-fam}
A family $\{ S_t\}_{t\in\T_+}$ of continuous  mappings of $X$ into itself
is said to be {\em evolution operator} (or evolution semigroup, or semiflow)
if  it satisfies the semigroup property:
\[
S_0=Id,~~~ S_{t+\tau}=S_t\circ S_\tau~~\mbox{for all}~~ t,\tau\ge 0.
\]
In the case when $\T=\R$ we assume in addition that the mapping $t\mapsto S_tx$ is continuous
from $\R_+$ into $X$ for every $x\in X$.  The pair $(X,S_t)$
is said to be a {\em dynamical system}
\index{dynamical system}
with the phase (or state)
space $X$ and the evolution operator $S_t$.
\par
If  $\T=\Z$, then  evolution operator (and dynamical system)  is
called discrete (or with discrete time).
If $\T=\R$, then $S_t$  (resp.\ $(X,S_t)$) is
called an  evolution operator
(resp.\ dynamical system) with continuous time.
\end{definition}

For any $D\subset X$ we denote by
$$
\gamma^t_D\equiv
\bigcup_{\tau\ge t}S_\tau D
$$
 the {\em tail} (from the moment $t$) of the  trajectories
\index{tail of trajectory} emanating from $D$. It is clear that
$\gamma^t_D=\gamma^0_{S_tD}\equiv\gamma^+_{S_tD}$. If $D=\{ v\}$ is
a single point set, then $\gamma^+_v\equiv \gamma^0_v$ is said to be
a   {\em positive semitrajectory (or semiorbit)} emanating from $v$.
 A curve  $\gamma \equiv \{ u(t)\,
:\, t\in\T\}$ in $X$ is said to be a {\em full trajectory}  iff
$S_tu(\tau)=u(t+\tau)$ for any $\tau\in\T$ and $t\ge 0$.

The set
\begin{equation}\label{7.1.3}
\omega(D)\equiv
\bigcap_{t>0}\overline{\gamma_D^t}
=
\bigcap_{t>0}\overline{\bigcup_{\tau\ge t}S_\tau D}
\end{equation}
is called the  {\em $\omega$-limit set} of the trajectories
emanating from $D$ (the bar over a set means the closure).

\begin{definition}[Global attractor]\label{de7.2.1}\index{global attractor}
\index{attractor!global}
Let $S_t$ be an evolution operator on a complete metric space $X$.
A bounded
 closed set $\fA\subset X$
is said to be a {\em global attractor} for $S_t$  if
\begin{enumerate}
\item[(i)] $\fA$ is an invariant set; that is,
 $S_t\fA =\fA$ for $t\ge 0$.
\item[(ii)] $\fA$ is uniformly attracting; that is,  for all bounded set $D\subset X$
\begin{equation}\label{7.2.1}
\lim_{t\to +\infty}
d_X\{S_tD\, |\, \fA\}=0\quad\mbox{for every bounded set $D\subset X$},
\end{equation}
where $d_X\{A | B\}= \sup_{x\in A} \dist_X (x, B)$ is the Hausdorff semidistance.
\index{Hausdorff semidistance (semimetric)}
\end{enumerate}
\end{definition}

\begin{definition}[Milnor's attractor]
Let a Borel measure $\mu$ such that $\mu(X)<\infty$ be given on the phase space $X$ of
a dynamical system $(X,S_t)$. 
A bounded closed set $\fA\subset X$
is said to be a  Milnor attractor
(with respect to the measure $\mu$ ) for $(X,S_t)$ if is a minimal closed invariant set
possessing the property
\begin{equation*}
  \lim_{t\arr\infty} \dist(S_t y, \fA)=0
\end{equation*}
for almost all (with respect to measure $\mu$) elements  $y\in X$.
\end{definition}

\section{Algorithms}

The construction of the cover of Milnor's the attractor is based on the following
approach. Next, we will describe an algorithm for a two-dimensional
domain, for a three-dimensional domain everything will be similar.

We consider a first order ODE on $\R^2$
\begin{equation} \label{gen_DS}
  \dot{x}=F(x), \qquad t\in (0,T)
\end{equation}
and suppose that it generates a dissipative dynamical system $(\R^2, S_t)$.

Suppose we have a rectangular region than contains absorbing set  

$\varOmega=\left[x_{min},x_{max}\right]\times\left[y_{min},y_{max}\right]$. 

Let's build on $\varOmega$ uniform rectangular
grid$\varOmega_{\text{h}}$:

$\varOmega_{h}=\left\{ \left(x_{i},y_{j}\right),x_{i}=x_{min}+ih_{x},\,y_{j}=y_{min}+jh_{y};i=0,...,N;j=0,...,M\right\} ,$

here $\text{N, M}$- number of points in the grid along the
axes $OX,\,OY$; $h_{x}=\frac{x_{max}-x_{min}}{N},h_{y}=\frac{y_{max}-y_{min}}{M}$

{\bf Computational algorithm}

\begin{enumerate}
\item We solve equation (?) taking a center of each cell $x_{ij}$ as initial state, for a time interval $(0,T)$ 
with a fixed time step $\Delta t$, using Runge-Kutta method. Thus, for each initial state we get a number of points $S_{k\cdot \Delta t}x_{ij}$, which represent the corresponding trajectory of the dynamical system.

\item In each subdomain $\varOmega_{\text{h}}$, we calculate a number of trajectory representatives,  that fall into this region

\item We filter the received data, and keep only those $\varOmega_{\text{h}}$, in which number of  trajectory representatives is bigger then a certain threshold value.

\item We divide every $\varOmega_{\text{h}}$ which we kept into $2^n$ ($n$ is a dimension of the phase space) smaller boxes.

\item We proceed form the step 1.
\end{enumerate}

{\bf Filtration}.

At present, we use the mean value filter:
$\varepsilon=\frac{1}{M*N}\sum_{i=1,j=1}^{N,M}I_{ij}$.

We discard  all the regions in which  a number of trajectories representatives is less then $\varepsilon$.

{\bf solution of ODE.} For a numerical solution of attractors, the Runge-Kutta method
of 4-order
\begin{align*}
y_{i+1} & =y_{i}+\frac{\varDelta t}{6}(k_{1}+2k_{2}+2k_{3}+k_{4})\\
 & k_{1}=f(x_{i},y_{i}),\,\\
 & k_{2}=f(x_{i}+0.5\varDelta t,y_{i}+0.5\varDelta tk_{1})\\
 & k_{3}=f(x_{i}+0.5\varDelta t,y_{i}+0.5\varDelta tk_{2})\\
 & k_{4}=f(x_{i}+\varDelta t,y_{i}+\varDelta tk_{3})
\end{align*}

The Runge-Kutta methods have a number of important advantages: 1)possess a sufficiently high degree of accuracy (with the
exception of the Euler method); 2)it is explicit, i.e. the value$y_{k+1}$ is calculated from the previously found values; 3) allow the use of a variable step, which makes it possible to reduce it where the function changes rapidly, and increase otherwise.

\section{Numerical experiments }

\begin{figure}[h!]
  \centering
  \includegraphics[width=\textwidth]{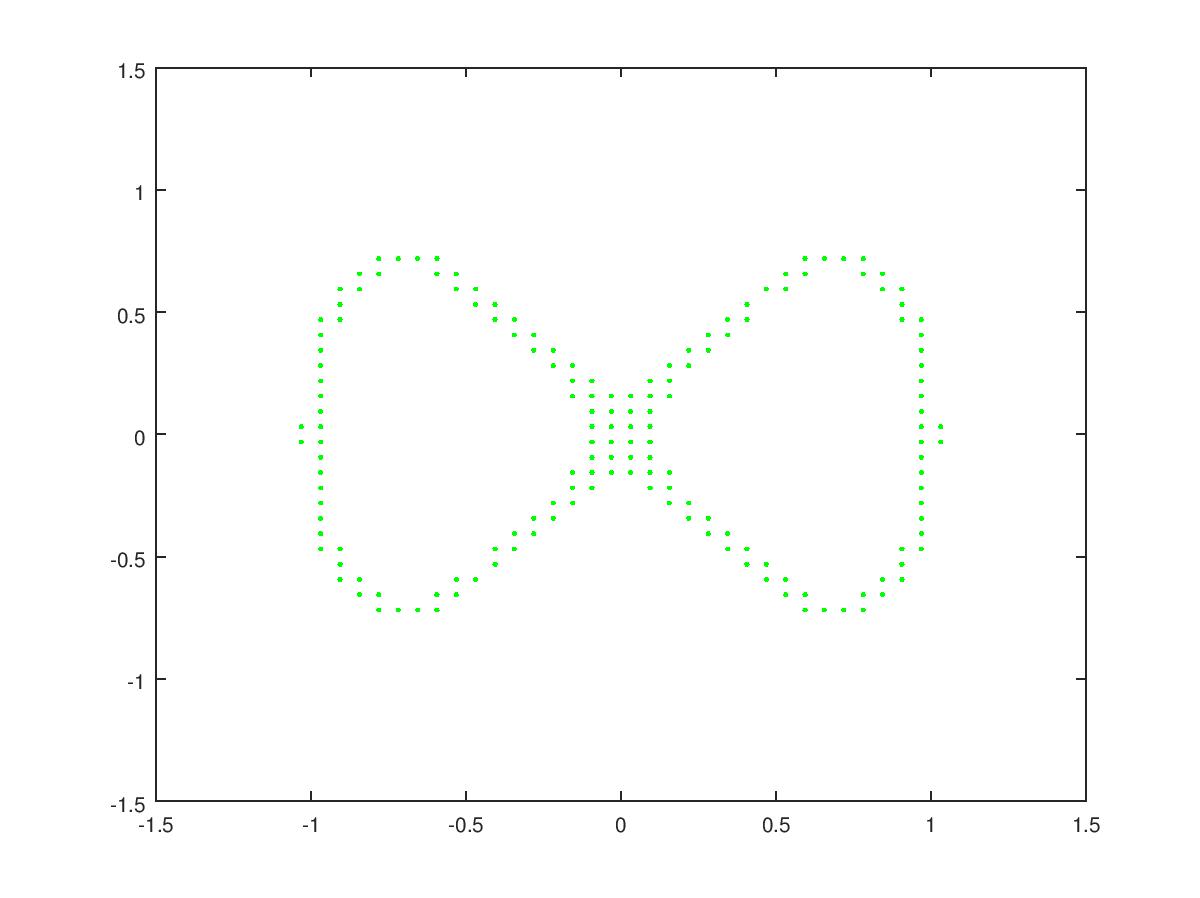}
  \caption{The first iteration of Algorithm 1 for the DS, generated by \eqref{eight_DS}}\label{eight_DS_1}
\end{figure}

\begin{figure}[h!]
  \centering
  \includegraphics[width=\textwidth]{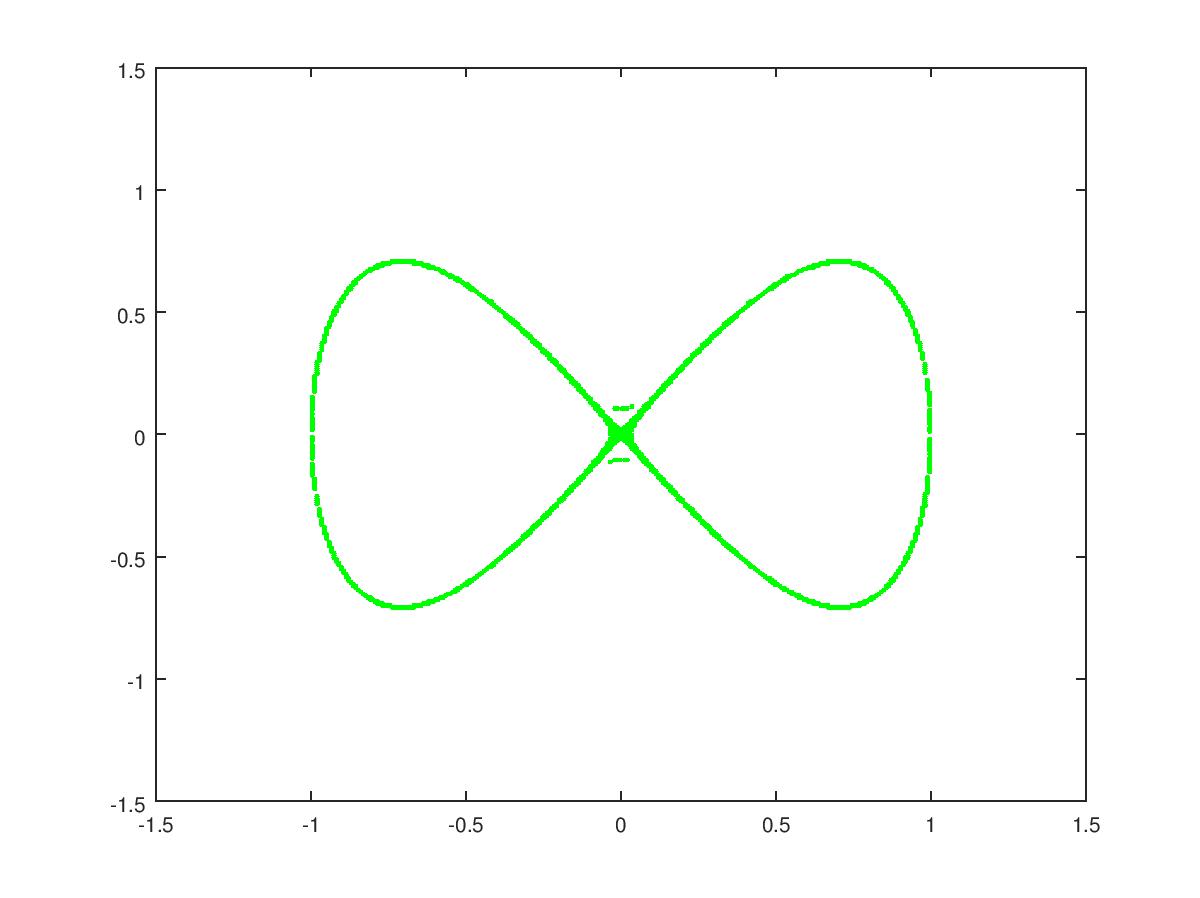}
  \caption{The forth iteration of Algorithm 1 for the DS, generated by \eqref{eight_DS}}\label{eight_DS_4}
\end{figure}

{\bf Example 1.} 2D dynamical system, generated by \eqref{eight_DS}. We take $x_{min}=-1.5,x_{max}=1.5,h_{x}=0.1$, $y_{min}=-1.5,y_{max}=1.5,h_{y}=0.1$, $\Delta t=0.01$ and time interval $(0,20)$. The first and the forth iterations are presented on figures \ref{eight_DS_1} and \ref{eight_DS_4}.

\begin{figure}[h!]
  \centering
  \includegraphics[width=\textwidth]{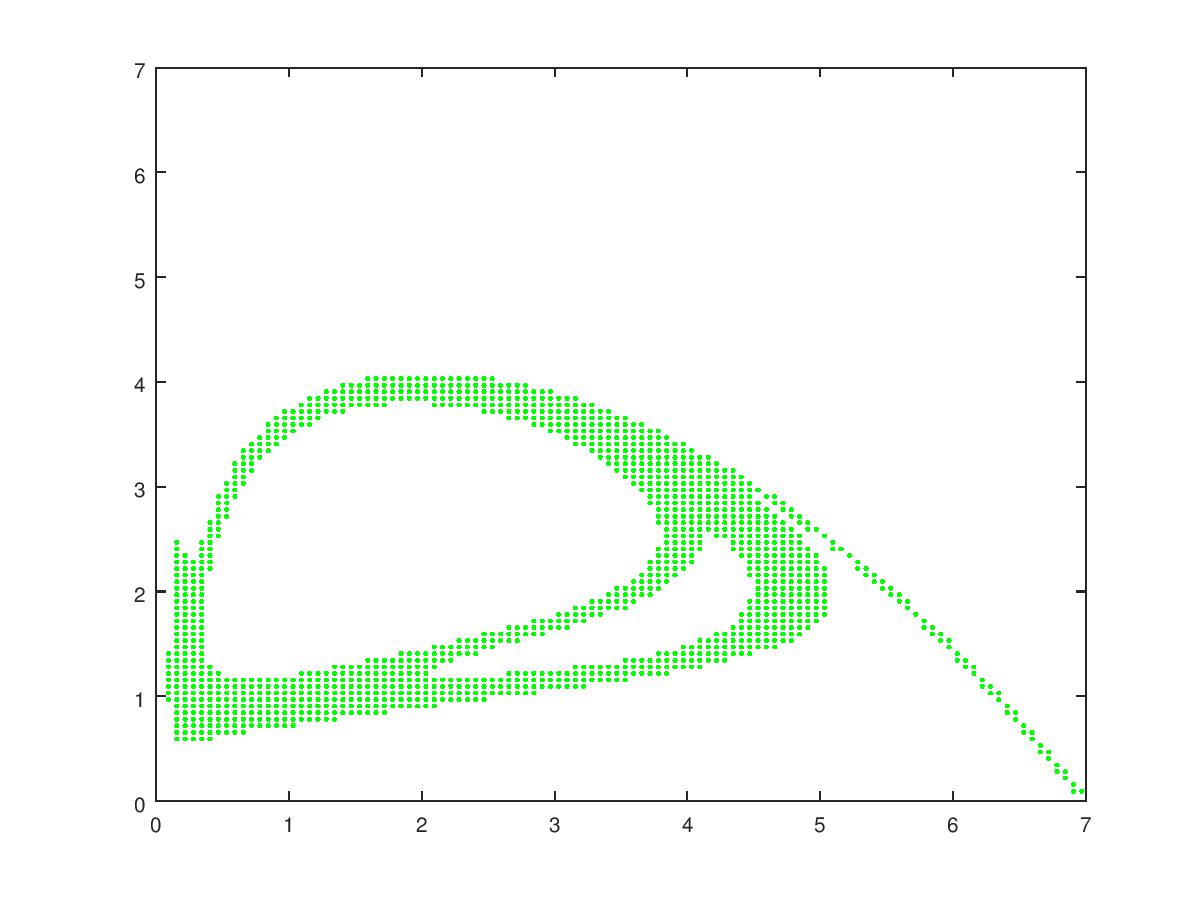}
  \caption{The first iteration of Algorithm 1 for the DS, generated by \eqref{HT}}\label{HT_1}
\end{figure}

\begin{figure}[h!]
  \centering
  \includegraphics[width=\textwidth]{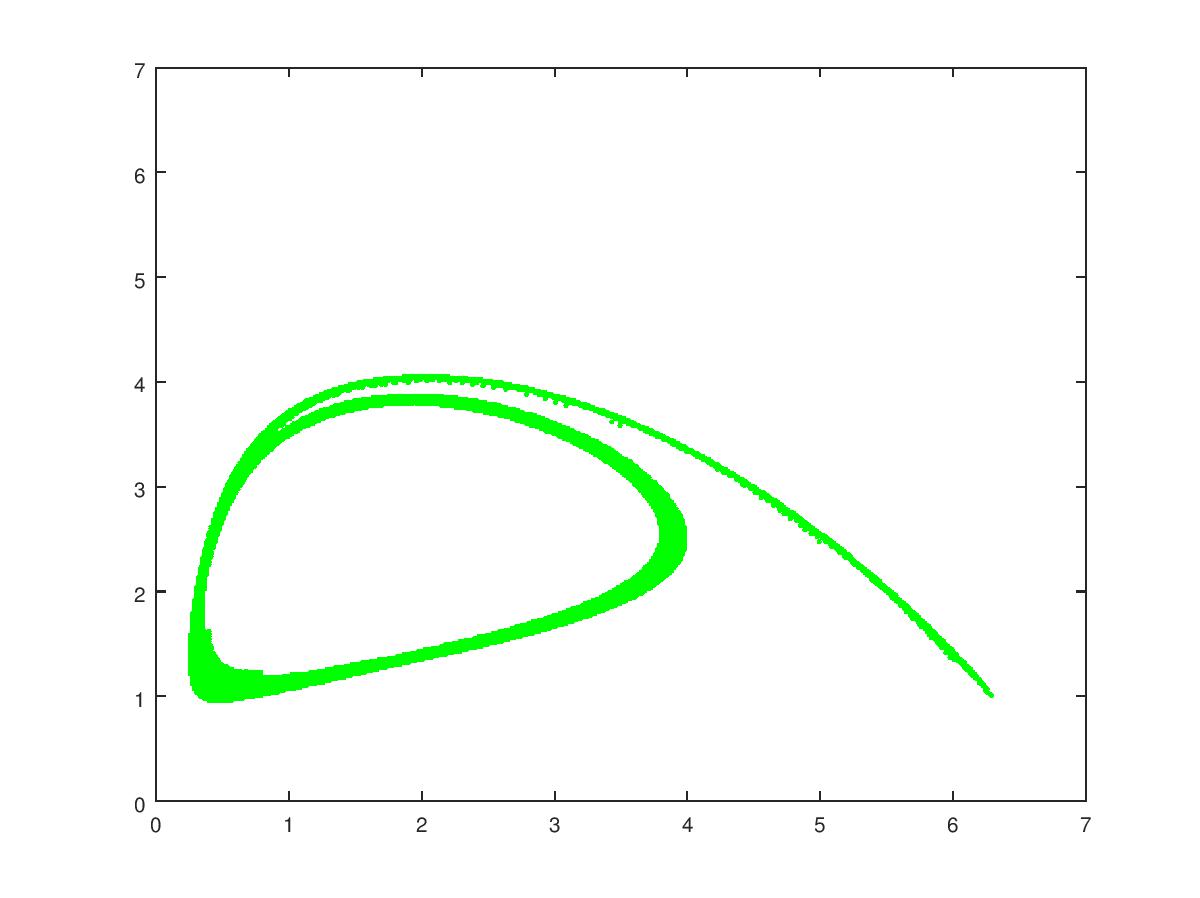}
  \caption{The third iteration of Algorithm 1 for the DS, generated by \eqref{HT}}\label{HT_3}
\end{figure}

\begin{figure}[h!]
  \centering
  \includegraphics[width=\textwidth]{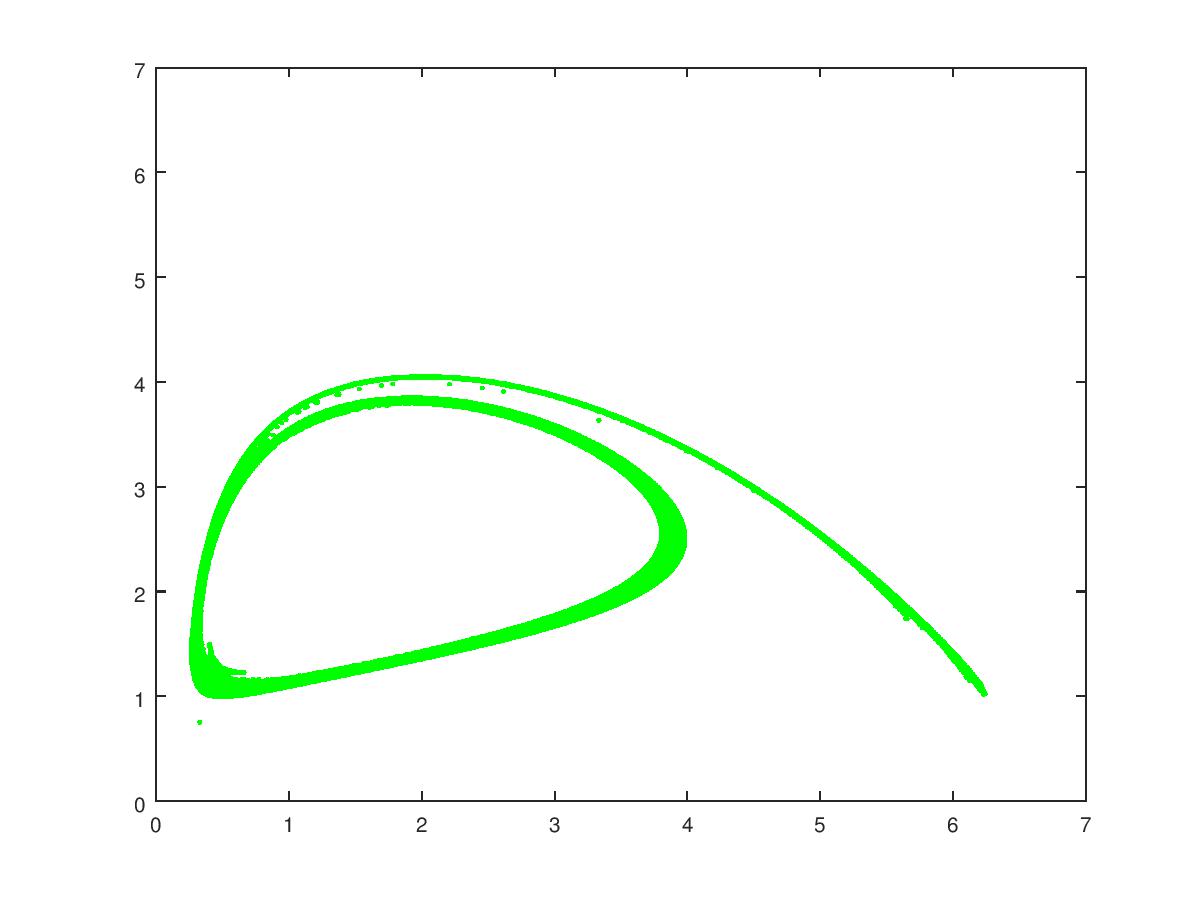}
  \caption{The fifth iteration of Algorithm 1 for the DS, generated by \eqref{HT}}\label{HT_5}
\end{figure}

\bigskip
{\bf Example 2.} Holling-Tanner model for predator-prey interaction \cite{Lyn2002}. We take $\Omega=(0,7) \times (0,7)$, $\Delta t=0.01$ and time interval $(0,40)$.
\begin{equation}\label{HT}
  \begin{aligned}
    \dot{x} =& x\left(1-\frac x7 \right) - \frac{6xy}{7+7x}\\
    \dot{y}= & 0.2y \left( 1-\frac{Ny}{x} \right)
  \end{aligned}
\end{equation}
We model the case $N=0.5$. The first, the  third and the fifth iterations are presented on figures \ref{HT_1}, \ref{HT_3} and \ref{HT_5}.

\bigskip

{\bf Example 3.} System \eqref{fp_appr}  is a 1-mode approximation of one fluid-structure interaction model. This is 3D system of ODE. We take $\Omega=(-3,3) \times (-1,1)\times (-1,1)$, $\Delta t=0.01$ and time interval $(0,10)$.

\begin{figure}[h!]
  \centering
  \includegraphics[width=\textwidth]{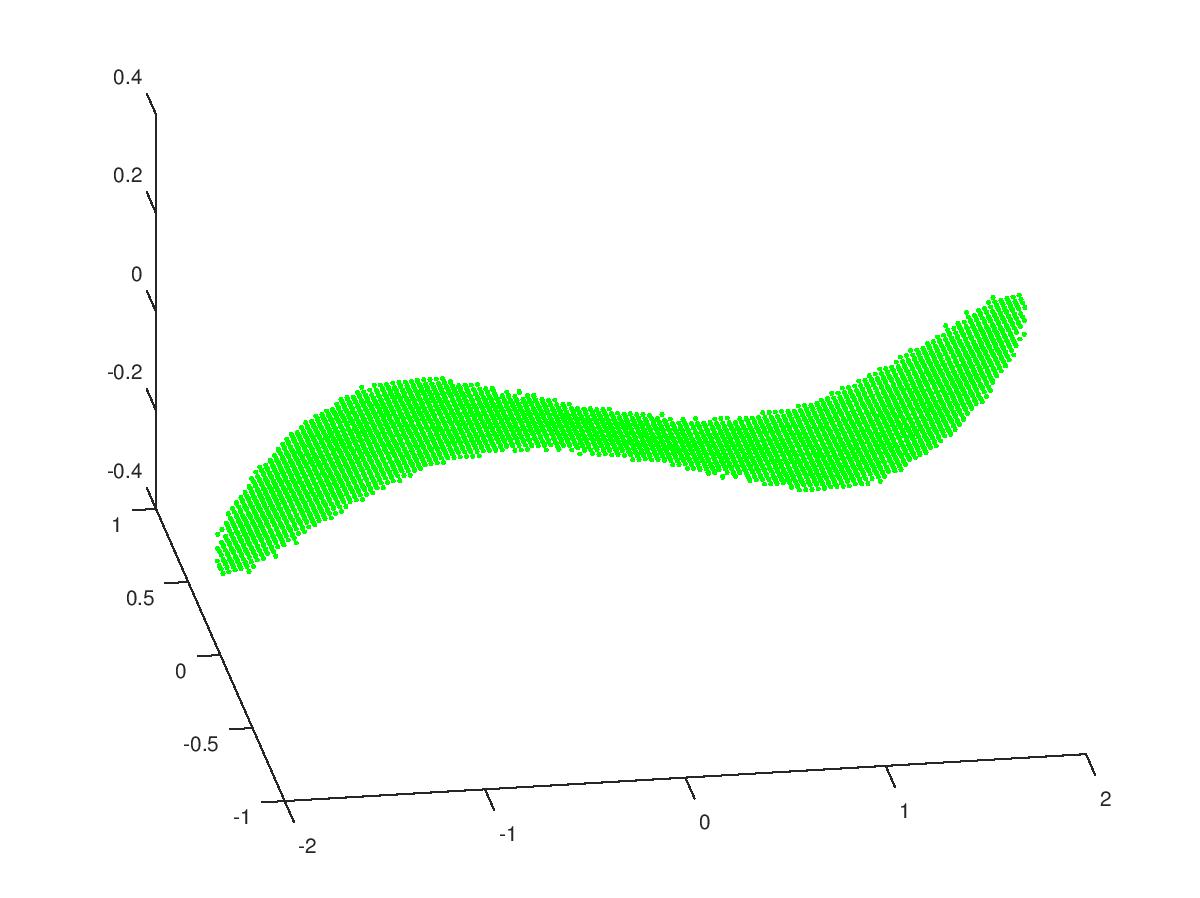}
  \caption{The second iteration of Algorithm 1 for the DS \eqref{fp_appr} }\label{FP_1}
\end{figure}

\begin{figure}[h!]
  \centering
  \includegraphics[width=\textwidth]{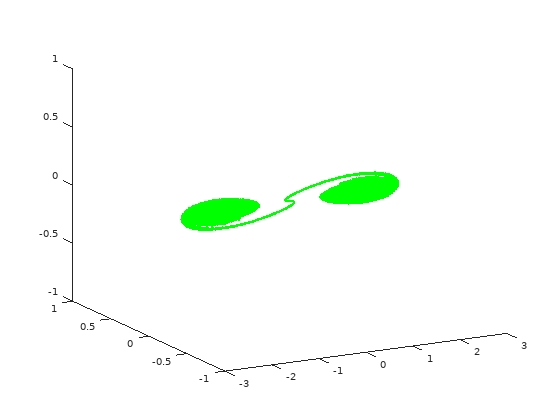}
  \caption{The ninth iteration of Algorithm 1 for the  DS \eqref{fp_appr} }\label{FP_9}
\end{figure}

{\bf Example 4.} Hopf's system \cite{Chu99}. This is 3D system of ODE. 
\begin{equation*}
  \begin{aligned}
    && \dot{x} + \mu x + y^2 + z^2  =  0,\\
    && \dot{y} + \nu y - xy - \beta z  =0, \\
    && \dot{z} +\nu z - xz +\beta y =0.
  \end{aligned}
\end{equation*}
where $\mu>0$, $\beta, \nu \in \R$. We perform numerical experiments for $\mu = 4$, $\beta=-1/4$, $\beta=1$. We take $\Omega=(-2,2) \times (-2,2)\times (-2,2)$, $\Delta t=0.01$ and time interval $(0,10)$.
\begin{figure}[h!]
  \centering
  \includegraphics[width=\textwidth]{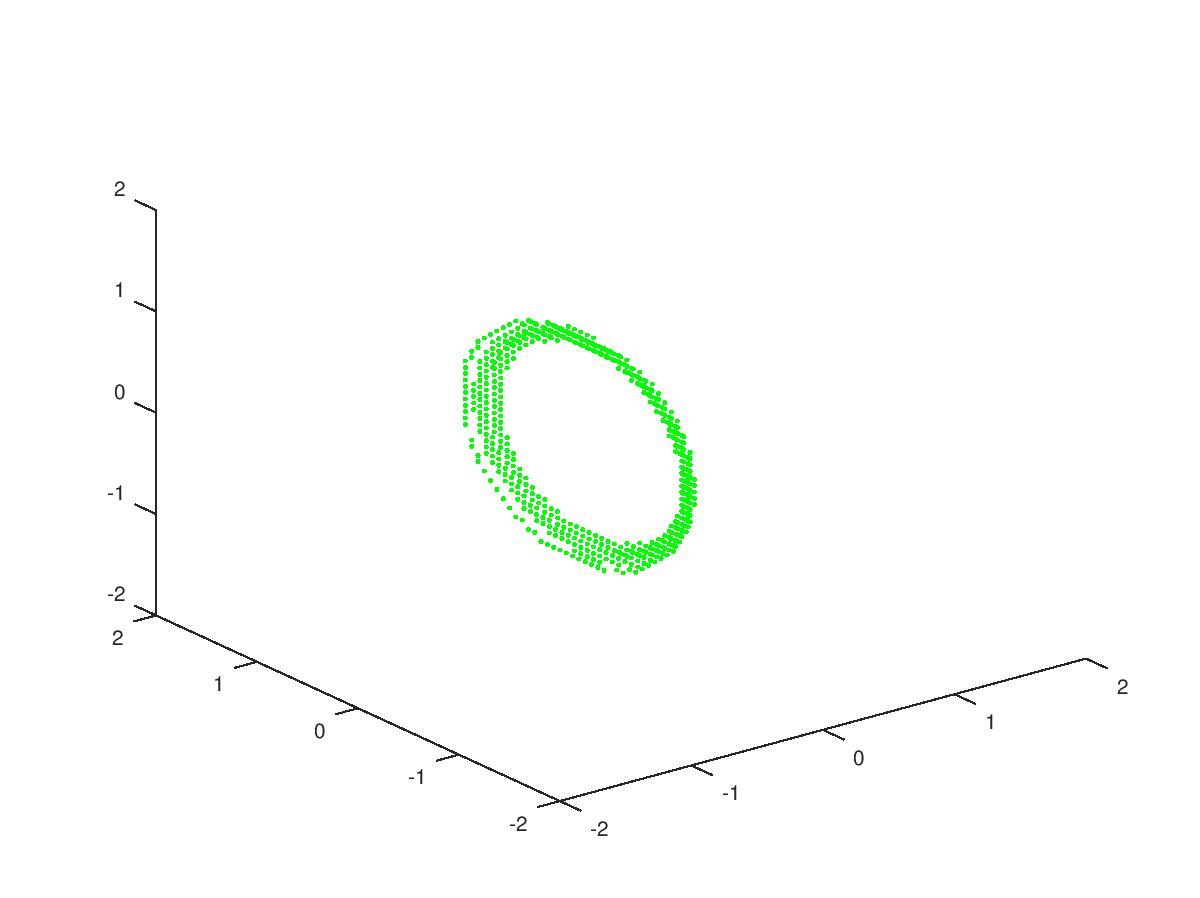}
  \caption{The first iteration of Algorithm 1 for the Hopf DS}\label{H_1}
\end{figure}

\begin{figure}[h!]
  \centering
  \includegraphics[width=\textwidth]{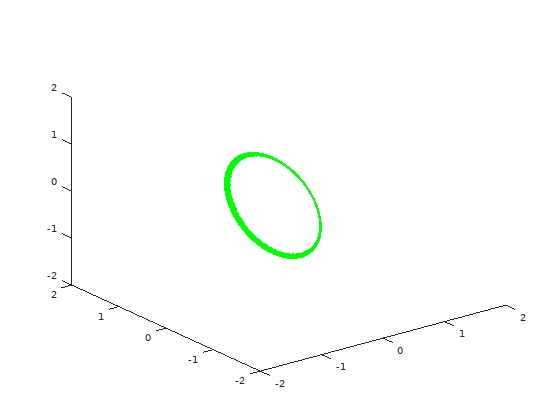}
  \caption{The fifth iteration of Algorithm 1 for the Hopf DS}\label{H_5}
\end{figure}

\section*{Acknowledgments} The authors were supported by the VolkswagenStiftung Project “Modeling, Analysis, and Approximation Theory toward Applications in Tomography and Inverse Problems”.

\end{document}